\documentclass[12pt]{article}
\newtheorem{theorem}{Theorem}
\newtheorem{lemma}{Lemma}
\font\tenmsbm=msbm10 scaled \magstep1
\font\sevenmsbm=msbm7  scaled \magstep1
\font\fivemsbm=msbm5  scaled \magstep1
\newfam\Bbbfam
\textfont\Bbbfam=\tenmsbm \scriptfont\Bbbfam=\sevenmsbm
\scriptscriptfont\Bbbfam=\fivemsbm
\def\Bbb{\fam\Bbbfam}

\begin{document}
\begin{center} 
Konstantin Igudesman \\
\vspace{0.5cm}
{\large \bf {STATISTICALLY SELF-SIMILAR FRACTAL SETS}}
\end{center}

In the present paper we define statistically self-similar sets, and, using a modification of method described in 
\cite{fal}, \cite{hut}, find a Hausdorff dimension of a statistically self-similar set.
 
\section{Hausdorff measure, dimention and metric} 

This section contains definitions and results we require (see \cite{fal} for more details).

If $U$ is a non-empty subset of ${\Bbb R}^n$, the diameter of $U$ is  
$$
|U|=\sup \{ |x-y| : x,y \in U \}.
$$

If $E\subset \bigcup \limits_{i} U_i $, and $0<|U_i|\leq \delta $ for each $i$, then $\{ U_i \}$ is a $\delta $-cover 
of $E$.  Let $E\subset {\Bbb R}^n$, and $s\geq 0$. For $\delta >0$, we define the measure
$$
{\cal H}^s_\delta (E)=\inf \sum \limits_{i=1}^{\infty }|U_i|^s , 
$$
where infimum is taken over all countable $\delta $-covers $\{ U_i \}$ of $E$.

The Hausdorff $s$-dimensional measure of $E$ is 
$$
{\cal H}^s (E)=\lim _{\delta \to 0}{\cal H}^s_\delta (E).
$$
The limit exists, but can be infinite, since ${\cal H}^s_\delta $ increases as $\delta $ decreases.

The $\sigma $-field of ${\cal H}^s$-measureable sets includes the Borel sets. Moreover, for any $E\subset 
{\Bbb R}^n$ there is a unique $s \in [0,n] $ such that
$$
{\cal H}^t (E)=
\left\{ 
\begin{array}{cc} 
\infty , & t<s  \\ 
0 ,& t>s   
\end{array} 
 \right. \ .
$$ 
The number $s$ is called the Hausdorff dimension (or fractal dimension) of $E$ and is denoted by $\dim (E)$.

If $E\subset {\Bbb R}^n$, and $\delta \geq 0$, the $\delta  $-parallel body of $E$ is the closed set 
$$
[E]_\delta =\{ x \in {\Bbb R}^n : \inf \limits _{y \in E} |x-y|\leq \delta  \}.
$$
Let $\Gamma $ be a collection of non-empty compact subset of ${\Bbb R}^n$.
The Hausdorff metric $d$ on $\Gamma $ is
$$
d(E,F)=\inf \{\delta : E\subset [F]_\delta , F\subset [E]_\delta \}.
$$
It is known that $\Gamma $ endowed with the Hausdorff metric $d$ is a complete metric space .

\section{Statistically self-similar sets} \label{stat}

In this section we construct a statistically self-similar set and find its Hausdorff dimension.
The collection of statistically self-similar sets contains the well-known self-similar sets.

A mapping $\psi : {\Bbb R}^n \to {\Bbb R}^n$ is called a similitude if $|\psi (x) - \psi (y)| = r |x-y|$ for all 
$x, y \in {\Bbb R}^n$, where $r<1$. Clearly, any similitude is a continuous mapping. 

Let $\{ \psi _j \}_{j=1}^m$ be a set of similitudes with ratios $\{ r_j \}_{j=1}^m$. 
We set $r=\max \limits _{1\leq j\leq m} r_j$. Let ${\cal N}_k \subseteq  \{ 1,\ldots ,m \}$ be a non-empty set,
for any positive integer $k$ (if ${\cal N}_k = \{ 1,\ldots ,m \}$ for all $k$, then the following construction gives a 
self-similar set). 
We set ${\cal L}^k=\prod \limits_{i=1}^{k} {\cal N}_i ,\  {\cal L}_p^k=\prod \limits_{i=p}^{k} {\cal N}_i \ (p\leq k)$. 
For any sequence $(j_1\ldots j_k) \in {\cal L}^k$ and any set $F$ we set 
$F_{j_1\ldots j_k} = (\psi _{j_1}\circ \ldots \circ \psi _{j_k})(F) $.
Now, let  $\Gamma $ be a collection of non-empty compact sets, and $\psi ^k : \Gamma \to \Gamma $, 
$$
\psi  ^k(F) = \bigcup \limits_{{\cal L}^k}  (\psi _{j_1}\circ \ldots \circ \psi _{j_k})(F) = 
\bigcup \limits_{{\cal L}^k} F_{j_1\ldots j_k}\ .
$$
for any positive integer $k$.

\begin{theorem} \label{conv}
There exists a unique set $E \in \Gamma $ such that for all $F \in \Gamma $
$$
\psi ^k(F) \stackrel{k \to \infty } \longrightarrow  E  
$$
with respect to the Hausdorff metric $d$. 
\end{theorem}

$\diamond $ 
If $F,G \in \Gamma $, then, by the definition  of $d$, 
$$
d(\psi ^k(F), \psi ^k(G))\leq \sup \limits _{{\cal L}^k} d(F_{j_1\ldots j_k}, G_{j_1\ldots j_k})  \leq r^k d(F,G) \ .
$$
Take any $F\in \Gamma $ and set $M=\sup \limits _{1\leq j \leq m} d(F,F_j)$. For any positive integers $p$ and 
$q$\  ($p<q$), 
$$
d(\psi ^p(F), \psi ^q(F))\leq 
d(\psi ^p(F), \psi ^p(\bigcup \limits_{{\cal L}_{p+1}^q} F_{j_{p+1}\ldots j_q})) \leq 
r^p d(F, \bigcup \limits_{{\cal L}_{p+1}^q} F_{j_{p+1}\ldots j_q}) \leq 
$$
$$
r^p \sup \limits_{{\cal L}_{p+1}^q} d(F, F_{j_{p+1}\ldots j_q}) \leq 
r^p \sup \limits_{{\cal L}_{p+1}^q} (d(F, F_{j_{p+1}})+d(F_{j_{p+1}}, F_{j_{p+1} 
j_{p+2}}) + \cdots +
$$
$$
d(F_{j_{p+1}\ldots j_{q-1}}, F_{j_{p+1}\ldots  j_{q}}) \leq 
r^p (M + r M + r^2 M + \cdots r^{q-p-1} M) < r^p \frac{M}{1-r}\ .
$$
Hence $d(\psi ^p(F), \psi ^q(F))\stackrel{p\to \infty }{\longrightarrow} 0$.
As $\Gamma $ endowed with the Hausdorff metric $d$ is a complete metric space, $\psi ^k(F)$ converges to 
a non-empty compact set $E$. 

Since for all $k$ and $G\in \Gamma$
$$
d(\psi ^k(G), E)\leq d(\psi ^k(G), \psi ^k(F))+d(\psi ^k(F), E)\ ,
$$
$E$ is unique.
$\diamond $

Let ${\cal K}$ be the set of finite sequences $\{ ( j_1 , \ldots , j_k )\}_{{\cal L}^k, k}$.
A finite subset ${\cal P}\subset {\cal K}$ is called a tree if for any $( j_1 , \ldots , j_k )\in {\cal K}$ there exists 
$( i_1 , \ldots , i_p )\in {\cal P}$ such that $j_q = i_q$ for all $q=1,\ldots , \min \{ k,p \}$, moreover, if $p\leq k$, 
then this $( i_1 , \ldots , i_p )\in {\cal P}$ is unique. Note that for any $k$ the set ${\cal L}^k$ is tree.

Given a tree ${\cal P}$, we set $p=\inf \{ k :  ( j_1 , \ldots , j_k ) \in {\cal P} \}$ (the length of the shortest branch 
of the tree) and $q=\sup \{ k :  ( j_1 , \ldots , j_k ) \in {\cal P} \}$ (the length of the longest branch of the tree).

\begin{lemma} \label{tree}
For any tree ${\cal P}$ and any non-negative numbers $\{ a_k \}_{k=1}^m$ we 
have
$$
\sum \limits_{{\cal P}} a_{j_1}\ldots a_{j_k} \geq 
\inf \limits_{p\leq k\leq q}\sum \limits_{{\cal L}^k} a_{j_1}\ldots a_{j_k}\ ,
$$
where $p$ is the length of the shortest branch of ${\cal P}$, and $q$ is the length of the longest one.
\end{lemma}

$\diamond $ 
For any $ ( j_1 , \ldots , j_p )\in {\cal L}^p$,\ 
${\cal P}_{j_1  \ldots  j_p} = \{ ( i_1 , \ldots , i_k )\in {\cal P} : ( j_1 , \ldots , j_p )=( i_1 , \ldots , i_p ) \}$.
If $( l_1 , \ldots , l_p )\in {\cal L}^p $, then 
\begin{equation} \label{tre}
{\cal P}' = \bigcup \limits_{{\cal L}^p}  
\bigcup \limits_{{\cal P}_{l_1  \ldots  l_p} } (j_1, \ldots , j_p , l_{p+1}, \ldots , l_k) 
\end{equation}
is a tree. Note that the first $p$ elements of
$(j_1, \ldots , j_p , l_{p+1}, \ldots , l_k)$ run through ${\cal L}^p$ and the remaining elements run through the 
ends of $( l_1 , \ldots , l_k ) \in {\cal P}_{l_1  \ldots  l_p}$.

For any $( j_1 , \ldots , j_p )\in {\cal L}^p$,
$$
\mu _{ j_1 \ldots  j_p} = 
\left\{ 
\begin{array}{ccc} 
\sum \limits_{{\cal P}_{j_1 \ldots  j_p}} a_{j_{p+1}}\ldots a_{j_k}  
& \mbox {if} & (j _1, \ldots , j_p) \not \in {\cal P} \\ 
1  & \mbox {if} &  (j _1, \ldots , j_p) \in {\cal P} \ .
\end{array} 
 \right. 
$$
Hence
\begin{equation} \label{mu} 
\sum \limits_{{\cal P}} a_{j_1}\ldots a_{j_k} = 
\sum \limits_{{\cal L}^p} a_{j_1}\ldots a_{j_p} \mu _{ j_1 \ldots  j_p}\ .
\end{equation}

Set $\mu = \inf  \limits_{{\cal L}^p} \mu _{ j_1 \ldots  j_p}$. From (\ref{mu}) it follows  that
$$
\sum \limits_{{\cal P}} a_{j_1}\ldots a_{j_k} \geq 
\sum \limits_{{\cal L}^p} a_{j_1}\ldots a_{j_k} 
$$
if  $\mu =1$.

If $\mu <1$, then there exists $(l _1, \ldots , l_p) \not \in {\cal P}$ such that 
$\mu = \mu _{ l_1 \ldots  l_p}  $. Thus, using (\ref{mu}), we get
$$
\sum \limits_{{\cal P}} a_{j_1}\ldots a_{j_k} \geq \mu _{ l_1 \ldots  l_p} \sum \limits_{{\cal L}^p} a_{j_1}\ldots 
a_{j_p} = \sum \limits_{{\cal L}^p} a_{j_1}\ldots a_{j_p} \sum \limits_{{\cal P}_{j_1 \ldots  j_p}} a_{l_{p+1}} \ldots 
a_{l_k} = 
$$
$$
\sum \limits_{{\cal P}'} a_{j_1} \ldots a_{j_p} a_{l_{p+1}} \ldots a_{l_k} \ .
$$
That follows from (\ref{tre}) that ${\cal P}'$ is a tree. Clearly, the length of the shortest branch of ${\cal P}'$ is 
equal to $p'>p$, and the length of the longest branch of ${\cal P}'$ is equal to $q' \leq q$. 

We proceed in this way until, after a finite number of steps, we reach a tree ${\cal L}^k$, where $p\leq k \leq q$.
$\diamond $ 

\begin{lemma} \label{ball}
Let $\{ V_i \}$ be a collection of disjoint open subset of ${\Bbb R}^n$ such that each $V_i $ contains a ball of 
radius $c_1 \rho $ and is contained in a ball of radius $c_2 \rho $. Then any ball $B$ of radius $\rho $ 
intersects, at most, $(1+2c_2)^n c_1^{-n}$ of the sets $\bar {V_i}$ (the bar denotes closure). 
\end{lemma}

$\diamond $ 
If $\bar {V_i}\bigcap B \neq \emptyset $,\  $\bar {V_i}$ is contained in a ball concentric with $B$ and of radius 
$(1+2c_2) \rho $. Let $h$ elements of the collection $\{\bar V_i \}$ intersect $B$, then summing up the 
volumes of the corresponding interior balls, we get $h (c_1 \rho )^n \leq (1+2c_2)^n \rho ^n$.
$\diamond $ 

\begin{lemma} \label{dim}
There exists a unique number $s\geq 0$ such that 
$$
\liminf \limits_{k\to \infty } \sum \limits_{{\cal L}^k} (r_{j_1}\ldots r_{j_k})^t= 
\left\{ 
\begin{array}{cc} 
\infty , & t<s  \\ 
0 ,& t>s   
\end{array} 
 \right. \ .
$$
\end{lemma}

$\diamond $ 
It is clear that $\liminf \limits_{k\to \infty } \sum \limits_{{\cal L}^k} (r_{j_1}\ldots r_{j_k})^t$ does not increase as 
$t$ increases from $0$ to $\infty $. Furthemore, if $t<t'$, then 
$$
\inf \limits_{l \geq k} \sum \limits_{{\cal L}^l} (r_{j_1}\ldots r_{j_l})^t =
\inf \limits_{l \geq k} \sum \limits_{{\cal L}^l} (r_{j_1}\ldots r_{j_l})^{t'+(t-t')} \geq 
(r^l)^{(t-t')} \inf \limits_{l \geq k} \sum \limits_{{\cal L}^l} (r_{j_1}\ldots r_{j_l})^{t'}\ , 
$$
which implies the statement.
$\diamond $ 

We say that the set of similitydes $\{ \psi _j \}_{j=1}^m$ satisfies the open set condition if there exists a 
bounded open set $V \subset {\Bbb R}^n$ such that 
\begin{equation} \label{open}
\psi _j(V)\subset V \  \mbox{and} \ \psi _j(V)\bigcap \psi _i(V)=\emptyset 
\end{equation}
for any $i, j=\overline {1,m}, \ j\neq i$. Thus the sets $\{ V_{j_1 \ldots j_k} \}_{{\cal L}^k, k}$ form a net  in the 
sence  that any two of sets  from the collection are either  disjoint, or one set is included into the other one.
The collection $\{ V_{j_1 \ldots j_k} \}_{\cal P}$ is disjoint for any tree $\cal P$.

If (\ref{open}) holds, then $\{\psi ^k(\bar V)\}_k$ is a decreasing sequence of compact sets, which convergs to 
$E$ with respect to the Hausdorff metric by Theorem \ref{conv}. It follows from the definition of Hausdorff 
metric that  $E=\bigcap \limits_{k=1}^{\infty }\psi ^k(\bar V)$.

The set $E$ from Theorem \ref{conv} is called statistically self-similar, if the open set conditions holds true.

\begin{theorem} \label{main}
The Hausdorff dimension of a statistically self-similar set $E$ is equal to $\dim (E) = s$ (see lemma \ref{dim}).
\end{theorem}

$\diamond $
For any ${\cal L}^k$ (moreover, for any tree ${\cal P}$) the collection $\{ \bar {V}_{j_1 \ldots j_k} \}_{{\cal L}^k}$ 
is a cover of $E$. 
As 
$$
\liminf \limits_{k\to \infty } \sum \limits_{{\cal L}^k} | \bar {V}_{j_1\ldots j_k} |^t = 
| \bar {V} | \liminf \limits_{k\to \infty } \sum \limits_{{\cal L}^k} (r_{j_1}\ldots r_{j_k})^t =
\left\{ 
\begin{array}{cc} 
\infty , & t<s  \\ 
0 ,& t>s   
\end{array} 
 \right. 
$$
and $| \bar {V}_{j_1 \ldots j_k} | \leq r^k |V| \stackrel{k \to \infty }{\longrightarrow} 0$, we get $\dim (E)\leq s$.

To prove the opposite inequality we show that, if $t \in [0,s)$, then ${\cal H}^t (E) = \infty $. Since $E$ is 
compact, it is sufficient to prove that
$$
\liminf \limits _{\delta \to 0} \sum | U_i |^t = \infty \ , 
$$
where the infimum is taken over all finite $\delta $-covers $\{ U_i \}$ of $E$. Given any $\delta $-cover 
$\{ U_i \}_{i=1}^N$ of $E$, we can cover $E$ by balls $\{ B_i \}_{i=1}^N$ with $| B_i |\leq 2 | U_i |$, then 
\begin{equation} \label{f1} 
\sum \limits_{i=1}^{N} | U_i |^t \geq  2^{-t} \sum \limits_{i=1}^{N} | B_i |^t \ .
\end{equation}

Suppose that an open set $V$ such that (\ref{open}) holds true contains a ball of radius $c_1$ and is 
contained in a ball of radius $c_2$. Take any $\rho \in ( 0 ,1 )$. For  each infinite sequence $(j_1, j_2, \ldots )$ 
with $(j_1,  \ldots ,j_k)\in {\cal L}^k$ for all $k$, curtail the sequence at the least value of $k$ such that
\begin{equation} \label{ro}
( \min \limits _{1\leq j \leq m} r_j ) \rho \leq r_{j_1} \ldots r_{j_k} \leq \rho \ ,
\end{equation}
and let us denote by $\cal P$ the set of finite sequences obtained in this way. It is clear that $\cal P$ is a tree.

Each $V_{j_1 \ldots j_k}$ contains a ball of radius $c_1 r_{j_1} \ldots r_{j_k}$ and hence a ball of radius 
$c_1\rho (\min \limits _{1\leq j \leq m} r_j ) $, by (\ref{ro}), and is contained in a ball of radius 
$c_2 r_{j_1} \ldots r_{j_k}$ and therefore of radius $c_2 \rho $. By lemma \ref{ball}, any ball $B$ of radius 
$\rho $ intersects, at most, $h = (c_1+2c_2)^n (c_1 \min \limits _{1\leq j \leq m} r_j )^{-n}$ sets of collection 
$\{ \bar {V}_{j_1 \ldots j_k} \}_{{\cal P}}$. Note that $h$ does not depend on $\rho $.

Let $B$ be a ball of radius $\rho $. We denote 
$$
{\cal D} = \{ ( j_1, \ldots , j_k) \in {\cal P}\  : \  \bar {V}_{j_1 \ldots j_k} \bigcap B \neq \emptyset  \}\ .
$$
The set $\cal D$ contains, at most, $h$ elements.
Hence
$$
\sum \limits_{\cal D} (r_{j_1} \ldots r_{j_k})^t \leq h \rho ^t = h | B |^t \ ,
$$
hence
\begin{equation} \label{f2} 
|B |^t \geq h^{-1} \sum \limits_{\cal D} (r_{j_1} \ldots r_{j_k})^t \ .
\end{equation}

Suppose that for any ball $B_i$ from (\ref{f1}) we have constructed a tree ${\cal P}_i$ and a subset 
${\cal D}_i \subset {\cal P}_i$ . It follows from (\ref{f2}) that
\begin{equation} \label{f3}
\sum \limits_{i=1}^{N} |B_i |^t \geq h^{-1} \sum \limits_{i=1}^{N} \sum \limits_{{\cal 
D}_i} (r_{j_1} \ldots r_{j_k})^t \ .
\end{equation}

Let ${\cal S} = \bigcup \limits_{i=1}^{N} {\cal D}_i $. From the construction of the $\cal S$ we see that the 
collection $\{ \bar {V}_{j_1\ldots j_k} \}_{\cal S}$ is a $(2c_2\delta )$-cover of $E$. We set
$$
{\cal P} = \{ ( j_1, \ldots , j_k)\in {\cal S} \ :\ \mbox{for any}\  l < k \  ( j_1, \ldots , j_l) 
\not \in {\cal S} \}\ .
$$
Clearly, ${\cal P}\subset {\cal S}$ is a tree, hence the collection $\{ \bar {V}_{j_1\ldots j_k} \}_{\cal P}$ is a 
$(2c_2\delta )$-cover of $E$.

It follows from lemma \ref{tree} that
\begin{equation} \label{f4} 
\sum \limits_{i=1}^{N} \sum \limits_{{\cal D}_i} (r_{j_1} \ldots r_{j_k})^t \geq 
\sum \limits_{\cal P} (r_{j_1} \ldots r_{j_k})^t \geq 
\inf \limits_{p \leq k \leq q} \sum \limits_{{\cal L}^k} (r_{j_1} \ldots r_{j_k})^t 
\end{equation}
where $p$ and $q$ are the lengths of the shortest and the longest branches of the tree ${\cal P}$, 
respectively.

Thus by (\ref{f1}), (\ref{f3}), (\ref{f4}), 
$$
\sum \limits_{i=1}^{N} | U_i |^t \geq 
2^{-t} h^{-1} \inf \limits_{p \leq k \leq q} \sum \limits_{{\cal L}^k} (r_{j_1} \ldots r_{j_k})^t \ .
$$

As $\rho _i \to 0$ as $\delta \to 0$, we have $p \to \infty $. 
It follows from the equality $\lim \limits_{k \to \infty }\sum \limits_{{\cal L}^k} (r_{j_1} \ldots r_{j_k})^t = \infty $,  
$t<s$, that ${\cal H}^t (E) = \infty $.
$\diamond $

\section {The intersection of Cantor's sets}

Let $K_0 = [0,1], \ K_1 = [0,\frac{1}{3}] \bigcup [\frac{2}{3} ,1],\ K_2 = [0,\frac{1}{9}] \bigcup 
[\frac{2}{9} ,\frac{1}{3}] \bigcup [\frac{2}{3} ,\frac{7}{9}] \bigcup [\frac{8}{9} ,1],$ etc., where $K_{i+1}$ is obtained 
by removing the open middle third of each interval in $K_i$. Then $K_i$ consists of $2^i$ intervals, each of 
length $3^{-i}$. Cantor's set is the perfect compact set $K = \bigcap \limits_{i=0}^{\infty } K_i $.

Consider the similitudes of the real line $\psi _1 (x) = \frac{x}{3},\  \psi _2 (x) = \frac{x}{3} + \frac{2}{3} $.
Since
$$
K_i = \bigcup \limits_{j_1\ldots j_i} (\psi _{j_1}\circ \ldots \circ \psi _{j_i})([0,1]) \ ,
$$
where the sum is taken over all $i$-tuples $\{ j_1\ldots j_i \}$, the Cantor set is the self-similar set (it is the 
particular case of a statistically self-similar set). The Hausdorff dimension of the Cantor set $K$ is equal to 
$\frac{\ln 2}{\ln 3}$. 

Given $a \in [0,1]$, we set $K+a=\{ x : x-a \in K\}, E_a = K\bigcap (K+a)$. 

\begin{theorem} \label{kan}
For almost all $a\in [0,1]$ with respect to the Lebesgue measure,  $E_a$ is a statistically 
self-similar set and its Hausdorff dimension is equal to $\frac{\ln 2}{3 \ln 3}$. 
\end{theorem}

$\diamond $
The set $S_{a,i} = K_i \bigcap (K_i + a)$ is a cover of $E$ for any non-negative integer $i$. It is clear that $\{ 
S_{a,i} \}$ is the decreasing sequence of the sets, and $E_a = \bigcap \limits_{i=0}^{\infty } S_{a,i}$. From 
definition of $\delta $-parallel body we conclude that $S_{a,i}\subset [S_{a,i+1}]_{3^{-i}}$. Hence
\begin{equation}\label{eq1}
d(S_{a,i}, E_a)\stackrel{i\to \infty }{\longrightarrow} 0\ .
\end{equation}

Let $a=0,a_1a_2a_3\ldots $ be the triadic expansion of $a \in [0,1]$. We can exclude from the consideration 
the set of $a$ for which the triadic expansion is ambiguously determined, because the Lebesgue measure of 
this set, is equal to zero (for these $a$ the set $E_a$ is either finite or similar to Cantor's set). For any 
$a \in [0,1]$ we dissect the set $\{a_i \}_{i=1}^{\infty }$ to subsets 
$A(a)=\{ a_k | \sum \limits_{i=0}^{k-1}a_i=0\  \rm{(mod\  2)}\}$, 
$\tilde A(a)=\{ a_k | \sum \limits_{i=0}^{k-1}a_i=1 \ \rm{(mod\  2)}\}$ where $a_0=0$.
Let ${\cal N}_{a,k}$ be a subset of $\{ 1,2 \} $ defined by 
$$
{\cal N}_{a,k} = \left\{
\begin{array}{cccc} 
\{ 1 \} & \mbox{if} & a_{k} \neq  2,\ a_{k}\in \tilde A(a), & \\ 
\{ 2 \} & \mbox{if} & a_{k} \neq  0,\ a_{k}\in A(a), & \\ 
\{ 1,2 \} & \mbox{if} & a_{k} = 2,\ a_{k}\in \tilde A(a)&\!\!\!\! \mbox{or}\ a_{k} = 0,\ a_{k}\in A(a)\ , 
\end{array}  
\right. 
$$
for any $a \in [0,1]$ and any positive integer $k$.

The similitudes $\psi _1 (x) = \frac{x}{3},\  \psi _2 (x) = \frac{x}{3} + \frac{2}{3} $ with ratios 
$r_1 = r_2 = \frac{1}{3}$ and sets ${\cal N}_{a,i}$ generate (see. Section \ref{stat}) a statistically self-similar set 
$G_a$ (the open set condition holds for the open interval $(-1, 2)$). As Section \ref{stat} we set
$$
{\cal L}^k_a=\prod \limits_{i=1}^{k} {\cal N}_{a,i} \ ,\quad 
\psi  ^k_a(F) = \bigcup \limits_{{\cal L}^k_a} (\psi _{j_1}\circ \ldots \circ \psi _{j_k})(F) = 
\bigcup \limits_{{\cal L}^k_a} F_{j_1\ldots j_k}\ ,
$$
where $F$ is a non-empty compact set. By Theorem \ref{conv}, 
\begin{equation}\label{eq2}
d(G_a, \psi  ^k_a([0,1]))\stackrel{k\to \infty }{\longrightarrow} 0\ .
\end{equation}

It follows from the construction that $S_{a,k}\subset \psi  ^k_a([0,1])$ and 
$\psi  ^k_a([0,1])\subset [S_{a,k}]_{3^{-k}}$, hence
\begin{equation}\label{eq3}
d(\psi  ^k_a([0,1]), S_{a,k})\stackrel{k\to \infty }{\longrightarrow} 0\ .
\end{equation}

Thus by (\ref{eq1}),\ (\ref{eq2}),\ (\ref{eq3}), 
$$
d(G_a, E_a) \leq d(G_a, \psi ^k_a ([0,1])) +d (\psi ^k_a ([0,1]), S_{a,k}) + d(S_{a,k} , E_a) \stackrel{k\to 
\infty }{\longrightarrow} 0\ ,
$$ 
so $E_a = G_a$.

To find the Hausdorff dimention of $E_a$ it is sufficient to calculate $s_a$ for which  
\begin{equation} \label{sum} 
\liminf \limits_{k\to \infty } \sum \limits_{{\cal L}^k_a} (r_{j_1}\ldots r_{j_k})^t  = 
\liminf \limits_{k\to \infty } 3^{-tk} M_a (k) = 
\left\{ 
\begin{array}{cc} 
\infty , & t<s_a  \\ 
0 ,& t>s_a   
\end{array} 
 \right.
\end{equation}
where $M_a (k)$ is the number of elements of ${\cal L}^k_a$.

Let $\chi _Z$ be the characteristic function of a set $Z$. If $Z$ is a non-negative integer, we consider it as a 
subset of the set of all non-negative integers. By the definition of ${\cal L}^k_a$,
\begin{equation}\label {part}
M_a(k)=2^{f_a(k)},\quad \mbox{where} \quad 
f_a(k)=\sum \limits_{i=1}^{k} (\chi _0(a_i)\chi _{A(a)}(a_i)+\chi _2(a_i)\chi _{\tilde A(a)}(a_i)) \ .  
\end{equation}

Clearly, if the triadic expansion of $a\in [0,1]$ is known, then we can calculate the Hausdorff dimension of 
$E_a$ by (\ref {sum}) and (\ref {part}), but for an arbitrary irrational number $x$ we cannot calculate the part of 
$0,\ 1$ or $2$ in the triadic expansion of $x$. The only known fact is that, for almost all $x\in [0,1]$ (with 
respect to the Lebesgue measure) the part of $0,\ 1$ and $2$ in the triadic expansion of $x$ is equal to 
$\frac{1}{3}$.

Let $\Omega =\{ \omega =(a_1,a_2,a_3,\ldots )\}$ be the set of sequences of independent random variates, 
where $a_i$ takes values $0,1,2$ with probability $\frac{1}{3}$. We define a mapping  
$a:\Omega \rightarrow [0,1], \  a=a(\omega )=\sum \limits_{k=1}^{\infty }\frac{a_k}{3^k}$,\  
then $\omega $ is the triadic expansion of $a$. The probability $P$ on $\Omega $ generates the measure 
$a_\#P$ on the $[0,1]$ by $a_\#P(B)=P(a^{-1}(B))$, where $B\subset [0,1]$. The measure $a_\#P$ coincides 
with the Lebesgue measure "mes" (see \cite{bor}).

Let $\gamma :\Omega \rightarrow \hat \Omega $ be the mapping defined by
$\gamma (\omega )=\hat \omega =(\hat a_1,\hat a_2,\hat a_3,\ldots )$, where 
$$
\begin{array}{ccc}
\hat a_k=0 & \mbox{if} & a_k=0,\  a_k\in A(a); \\
\hat a_k=1 & \mbox{if} & a_k=1,\  a_k\in A(a); \\
\hat a_k=2 & \mbox{if} & a_k=2,\  a_k\in A(a); \\
\hat a_k=3 & \mbox{if} & a_k=0,\  a_k\in \tilde A(a); \\
\hat a_k=4 & \mbox{if} & a_k=1,\  a_k\in \tilde A(a); \\
\hat a_k=5 & \mbox{if} & a_k=2,\  a_k\in \tilde A(a)\ .
\end{array}
$$
The probability $P$ on $\Omega $ generates the probability $\hat P=\gamma _\#P$ on $\hat \Omega $. The 
sequence of random variates $\hat a_1,\hat a_2,\hat a_3,\ldots $ forms a homogeneous Marcov chain with 
the starting distribution $p=(\frac{1}{3},\frac{1}{3},\frac{1}{3},0,0,0)$ and the matrix of transition probabilities 
$$
\| p_{i,j} \|=\left(
\begin{array}{cccccc}
\frac{1}{3} & \frac{1}{3} & \frac{1}{3} & 0 & 0 & 0 \\
0 & 0 & 0 & \frac{1}{3} & \frac{1}{3} & \frac{1}{3} \\
\frac{1}{3} & \frac{1}{3} & \frac{1}{3} & 0 & 0 & 0 \\
0 & 0 & 0 & \frac{1}{3} & \frac{1}{3} & \frac{1}{3} \\
\frac{1}{3} & \frac{1}{3} & \frac{1}{3} & 0 & 0 & 0 \\
0 & 0 & 0 & \frac{1}{3} & \frac{1}{3} & \frac{1}{3}
\end{array}
\right)_.
$$

It follows from the strong law of large numbers that 
$$
P \left\{ \omega  : \lim \limits_{k\rightarrow \infty }\frac{1}{k} 
\sum \limits_{i=1}^{k}\chi _0(a_i)=\frac{1}{3} \right\} = 1\ . 
$$
Hence
\begin{equation} \label{ver}
\hat P \left\{ \hat \omega  : \lim \limits_{k\rightarrow \infty }\frac{1}{k} 
\sum \limits_{i=1}^{k} (\chi _0(\hat a_i)+\chi _3(\hat a_i)) = \frac{1}{3} \right\} = 1\ . 
\end{equation}
Since the starting distribution and the matrix of transition probability do not fluctuate if we interchange the  
random variates $a_k = 3$ and $a_k = 5$, and the set in the braces in (\ref{ver}) depends only on 
the starting distribution and the matrix of transition probability, we have
$$
\hat P \left\{ \hat \omega  : \lim \limits_{k\rightarrow \infty }\frac{1}{k} 
\sum \limits_{i=1}^{k}(\chi _0(\hat a_i)+\chi _5(\hat a_i)) = \frac{1}{3} \right\} = 1\ . 
$$

It follows from $\chi _0(\hat a_i) = \chi _0(a_i) \chi _{A(a)}(a_i)$ and 
$\chi _5(\hat a_i) = \chi _2(a_i) \chi _{\tilde A(a)}(a_i)$ that 
$$
\rm {mes}\left\{ {a\in[0,1] : \lim \limits_{k\rightarrow \infty 
}\frac{f_a(k)}{k}=\frac{1}{3}}\right\} = 
\phantom{222222222222222222222222222222222222222}
$$
$$
\phantom{2222222}
a_\# P \left\{ {a\in[0,1] : \lim \limits_{k\rightarrow \infty }\frac{1}{k} 
\sum \limits_{i=1}^{k} (\chi _0(a_i)\chi _{A(a)}(a_i)+\chi _2(a_i)\chi _{\hat 
A(a)}(a_i))= \frac{1}{3}}\right\} = 
$$
$$
\phantom{2222222222222221}
P \left\{ {\omega  : \lim \limits_{k\rightarrow \infty }\frac{1}{k} 
\sum \limits_{i=1}^{k} (\chi _0(a_i)\chi _{A(a)}(a_i)+\chi _2(a_i)\chi _{\hat 
A(a)}(a_i))= \frac{1}{3}}\right\} = 
$$
$$
\phantom{2222222222222222222222222222}
\hat P \left\{ \hat \omega  : \lim \limits_{k\rightarrow \infty }\frac{1}{k} 
\sum \limits_{i=1}^{k}(\chi _0(\hat a_i)+\chi _5(\hat a_i)) = \frac{1}{3} \right\} = 1\ . 
$$
Hence, using (\ref{sum}), except for a set of $a\in [0,1]$ of the Lebesgue measure $0$, 
$$
\liminf \limits_{k\to \infty } \sum \limits_{{\cal L}^k_a} (r_{j_1}\ldots r_{j_k})^t  = 
\liminf \limits_{k\to \infty } (2^\frac{f_a(k)}{k} 3^{-t})^k = 
\left\{ 
\begin{array}{cc} 
\infty , & t < s_a = \frac{\ln 2}{3 \ln 3}  \\
\phantom{.} & \phantom{.} \\ 
0 ,& t > s_a = \frac{\ln 2}{3 \ln 3}\ _.    
\end{array} 
 \right. 
$$

It follows that for almost all $a\in [0,1]$,
$$
\dim (E_a) = \frac{\ln 2}{3 \ln 3}\ .\quad \diamond 
$$

It follows from Theorem \ref{kan} that Mandelbrot's conjecture about the intersection of fractal sets does not 
hold for the intersection of Cantor sets. According to the conjecture the codimention of the intersection must 
be almost sure equal to the sum of the codimentions, i.e. be equal to $(\frac{\ln 4}{\ln 3} - 1)$ in our case.

Supported by grant of RFFI \quad N 00-01-00308.

\end{document}